\documentclass[a4paper,11pt,leqno]{amsart}

\usepackage{a4wide}
\usepackage{enumerate,xspace}
\usepackage{amsmath,amssymb,hyperref,delarray,srcltx}
\usepackage[english]{babel}
\usepackage[english,linesnumbered]{algorithm2e}
\usepackage[utf8]{inputenc}
\usepackage{graphics}
\SetKwFunction{BDyDF}{BDyDF}
\SetKwFunction{Bound}{Bound}
\SetKwFunction{NDelta}{NDelta}
\SetKwFunction{OptimalT}{OptimalT}
\SetKwFunction{floor}{floor}
\SetKwFunction{GRHcheck}{GRHcheck}
\expandafter\ifx\csname SetAlgoLined\endcsname\relax\def\SetAlgoLined{\relax}\fi
\input{macros}
\input{abbrev}

\title[An improvement to an algorithm of
Belabas, Diaz y Diaz and Friedman]{An improvement to an algorithm of\\
Belabas, Diaz y Diaz and Friedman}

\author[L.~Greni\'{e}]{Lo\"{\i}c Greni\'{e}}
\address[L.~Greni\'{e}]{Dipartimento di Ingegneria gestionale,
    dell'informazione e della produzione\\
    Universit\`{a} di Bergamo\\
    viale Marconi 5\\
    24044 Dalmine (BG)
    Italy}
\email{loic.grenie@gmail.com}

\author[G.~Molteni]{Giuseppe Molteni}
\address[G.~Molteni]{Dipartimento di Matematica\\
    Universit\`{a} di Milano\\
    via Saldini 50\\
    20133 Milano\\
    Italy}
\email{giuseppe.molteni1@unimi.it}

\begin{document}
\bibliographystyle{amsalpha}
\begin{abstract}
In \cite{small-generators} Belabas, Diaz y Diaz and Friedman show a way to
determine, assuming the Generalized Riemann Hypothesis, a set of prime ideals
that generate the class group of a number field. Their method is efficient
because it produces a set of ideals that is smaller than earlier proved
results. Here we show how to use their main result to algorithmically produce
a bound that is lower than the one they prove.
\end{abstract}
\maketitle
\section{Introduction}
We refer the reader to the paper \cite{small-generators} for an outline of
Buchmann's algorithm.

Let \K be a number field of degree \nK, with $r_1$ (resp. $r_2$) real (resp.
pair of complex) embeddings. We denote \DK the absolute value of its
discriminant.
\begin{defi}
Let \WC be the set of functions $F\colon[0,{+\infty})\to\RM$ such that
\begin{itemize}
\item $F$ is continuous;
\item $\exists\varepsilon>0$ such that the function
$F(x)e^{(\frac12+\varepsilon)x}$ is integrable and of bounded variation;
\item $F(0)>0$;
\item $(F(0)-F(x))/x$ is of bounded variation.
\end{itemize}
Let then, for $T>1$, $\WC(T)$ be the subset of $\WC$ such that
\begin{itemize}
\item $F$ has support in $[0,\log T]$;
\item the Fourier cosine transform of $F$ is non-negative.
\end{itemize}
\end{defi}

The main result of \cite{small-generators} is, up to a minor reformulation:
\begin{theo}[\textbf{Belabas, Diaz y Diaz, Friedman}]\label{theoKB}
Let \K be a number field satisfying the Riemann Hypothesis for all
\L-functions attached to non-trivial characters of its ideal class group
$\Cl_\K$, and suppose there exists, for some $T>1$, an $F\in\WC(T)$ with
$F(0)=1$ and such that
\begin{multline}\label{theoeq}
2\sum_\pG\log\N\pG\sum_{m=1}^{+\infty}\frac{F(m\log\N p)}{\N\pG^{m/2}}
    >
        \log\Delta_\K-\nK\g- \nK\log(8\pi)-\frac{r_1\pi}2\\
        + r_1\int_0^{+\infty}\frac{1-F(x)}{2\cosh(x/2)}\d x+\nK\int_0^{+\infty}\frac{1-F(x)}{2\sinh(x/2)}\d x\ .
\end{multline}
Then the ideal class group of \K is generated by the prime ideals of \K
having norm less than $T$.
\end{theo}

The authors apply the result to the function $\frac1LC_L\ast C_L$ where
$L=\log T$, $\ast$ is the convolution operator and $C_L$ is the characteristic
function of $({-\frac L2},\frac L2)$, to get the
\begin{coro}[\textbf{Belabas, Diaz y Diaz, Friedman}]\label{coroKB}
Suppose \K is a number field satisfying the Riemann Hypothesis for all
\L-functions attached to non-trivial characters of its ideal class group
$\Cl_\K$, and for some $T>1$ we have
\begin{multline}\label{coroeq}
2\sum_{\substack{\pG,m\\\N\pG^m<T}}\frac{\log\N\pG}{\N\pG^{m/2}}\left(1-\frac{\log\N\pG^m}{\log T}\right)
    >
        \log\Delta_\K-\nK\left(\g+\log(8\pi)-\frac{c_1}{\log T}\right)\\
        - r_1\left(\frac\pi2-\frac{c_2}{\log T}\right)\ ,
\end{multline}
where
$$c_1=\frac{\pi^2}2\ ,\quad c_2=4C\ .$$
(Here $C=\sum_{k>=0}(-1)^k(2k+1)^{-2}=0.915965\cdots$ is Catalan's constant.)

Then the ideal class group of \K is generated by the prime ideals of \K
having norm less than $T$.
\end{coro}
Our aim is to find a good $T$ for the number field \K as fast as possible
exploiting the bilinearity of the convolution product.
\section{Setup}
We use the following definition to simplify a little bit the language.
\begin{defi}
A \emph{bound} for \K is an $L=\log T$ with $T$ as in Theorem~\ref{theoKB}.
\end{defi}
\subsection{Rewriting the theorem}
We begin by homogenizing Equation~\eqref{theoeq} and relaxing the requirement
$F(0)=1$ to $F(0)>0$ so that now the condition on the function is
\begin{multline}\label{homoeq}
2\sum_\pG\log\N\pG\sum_{m=1}^{+\infty}\frac{F(m\log\N p)}{\N\pG^{m/2}}
    >
        F(0)\left(\log\Delta_\K-\nK\g- \nK\log(8\pi)-\frac{r_1\pi}2\right)\\
        + r_1\int_0^{+\infty}\frac{F(0)-F(x)}{2\cosh(x/2)}\d x+\nK\int_0^{+\infty}\frac{F(0)-F(x)}{2\sinh(x/2)}\d x\ .
\end{multline}
\begin{defi}
Let \SC be the real vector space of even and compactly supported step functions
and, for $T>1$, let $\SC(T)$ be the subspace of \SC of functions supported in
$\left[{-\frac{\log T}2},\frac{\log T}2\right]$.
\end{defi}
\begin{defi}
For any integer $N>=1$ and positive real $\delta$ we define the subspace
$\SC(N,\delta)$ of $\SC(e^{2N\delta})$ made of functions which are constant
$\forall k\in\NM$ on $[k\delta,(k+1)\delta)$.
\end{defi}
The elements of $\SC(N,\delta)$ are thus step functions with fixed step width
$\delta$. If $N>=1$, $\delta>0$ and $T=e^{2N\delta}$ we have
\begin{subequations}
\begin{alignat}{1}
\SC(N,\delta)               & \subset\SC(T)\subset\SC   \\
\forall\Phi\in\SC(T),\ \ \frac1{\|\Phi\|_2^{\,2}}\Phi\ast\Phi
                    & \in\WC(T)\\
\SC(N,\delta)               &\subset\SC(N+1,\delta)\label{incnext}\\
\forall k>=1,\ \ \SC\left(kN,\frac\delta k\right)
                    & \subseteq\SC(N,\delta)\label{incmul}\ .
\end{alignat}
\end{subequations}

If, for some $T>1$, $\Phi\in\SC(T)$ and $F=\Phi\ast\Phi$ satisfies
\eqref{homoeq} then, according to Theorem \ref{theoKB}, $\Cl_\K$ is generated
by prime ideals \pG such that $\N\pG<T$. This leads us to define the linear
form $\ell_\K$ on $\SC\ast\SC$ by
\begin{multline*}
\ell_\K(F)=-2\sum_\pG\log\N\pG\sum_{m=1}^{+\infty}\frac{F(m\log\N p)}{\N\pG^{m/2}}
    +
        F(0)\left(\log\Delta_\K-\nK\g- \nK\log(8\pi)-\frac{r_1\pi}2\right)\\
        + r_1\int_0^{+\infty}\frac{F(0)-F(x)}{2\cosh(x/2)}\d x+\nK\int_0^{+\infty}\frac{F(0)-F(x)}{2\sinh(x/2)}\d x
\end{multline*}
and the quadratic form $q_\K$ on \SC by $q_\K(\Phi)=\ell_\K(\Phi\ast\Phi)$. We
can at this point give a weaker version of Theorem~\ref{theoKB} as
\begin{coro}\label{coroKBL}
Let \K be a number field satisfying \GRH and $T>1$. If the restriction of
$q_\K$ to $\SC(T)$ has a negative eigenvalue then $\Cl_\K$ is generated by
prime ideals \pG such that $\N\pG<T$.
\end{coro}
Note that $q_\K$ is a continuous function as a function from
$(\SC(T),\|.\|_1)$ to \RM. Therefore if $\log T$ is a bound for \K then there
exists an $L'<\log T$ such that $L'$ is a bound for \K. Note also that, in
terms of $T$, only the norms of prime ideals are relevant, which means that we
do not need the smallest possible $T$ to get the best result.
\begin{rem*}
If $T>1$ and $\Phi\in\SC(T)$, then for any $\varepsilon>0$ there exists $N>=1$,
$\delta>0$ and $\Phi_\delta\in\SC(N,\delta)$ such that
$\|\Phi\ast\Phi-\Phi_\delta\ast\Phi_\delta\|_\infty<=\varepsilon$ and
$e^{2N\delta}<=T$. Hence we do not loose anything in terms of bounds if we
consider only the subspaces of the form $\SC(N,\delta)$.
\end{rem*}
\subsection{Computing the integrals}
Let $T>1$ be a real, $L=\log T$ and $F_L=C_L\ast C_L$ where, as above, $C_L$ is
the characteristic function of $\left[-\frac L2,\frac L2\right]$. We readily
see that $F_L(x)=(L-x)C_{2L}(x)$ for any $x>=0$. We easily compute
\begin{align*}
\int_0^{+\infty}\frac{F_L(0)-F_L(x)}{2\cosh(x/2)}\d x
&=
    4C-4\Imm\dilog\left(\frac i{\sqrt T}\right)
\intertext{and}
\int_0^{+\infty}\frac{F_L(0)-F_L(x)}{2\sinh(x/2)}\d x
&=
    \frac{\pi^2}2-4\dilog\left(\frac1{\sqrt T}\right)+\dilog\left(\frac1T\right)
\end{align*}
where $C$ is Catalan's constant and $\dilog(x)$ is the dilogarithm function
normalized to be the primitive of $-\frac{\log(1-x)}x$ such that $\dilog(0)=0$
(this is the normalization of \cite{PARI}).
\subsection{A remark on the restriction of quadratic forms}\label{rkqf}
Let $q$ be a quadratic form on an $n$-dimensional vector space $V$ of
signature $(z,p,m)$. We can interpret $p$ (resp. $m$) as the dimension of a
maximal subspace on which $q$ is positive (resp. negative) definite while the
kernel of $q$ has dimension $z=n-p-m$.

Let $H$ be an hyperplane of $V$ and $q'$ the restriction of $q$ to $H$. A
maximal subspace on which $q'$ is definite is a subspace on which $q$ is
definite, thus the intersection of a maximal subspace on which $q$ is
definite with $H$. This means the signature $(z',p',m')$ of $q'$ will be such
that $p'<=p<=p'+1$ and $m'<=m<=m'+1$. Cases $p=p'+1$, $m=m'+1$ and $p=p'$,
$m=m'$ are both possible with $z=n-p-m=z'-1$ and $z=z'+1$ respectively.
\section{Improving the result}
\subsection{Basic bound}
We restate \cite[Section 3, p. 1191]{small-generators} which determines an
optimal bound for Corollary~\ref{coroKB}. Let~$\GRHcheck(\K,\log T)$ be the
function that returns the right hand side of \eqref{coroeq} minus its left hand
side and $\BDyDF(\K)$ be the function which computes the optimal bound, by
dichotomy for instance. The computation of $\BDyDF(\K)$ is very fast because
the only arithmetic information we need on $\K\simeq\quot{\QM[x]}{(P)}$ is the
splitting information for primes $p<T$ and is determined easily for nearly all
$p$. Indeed if $p$ does not divide the index of $\quot{\ZM[x]}{(P)}$ in \OK,
then the splitting of $p$ in \K is determined by the factorization of $P\mod
p$. We can also store such splitting information for all $p$ that we consider
and do not recompute it each time we test whether a given bound $\log T$ is
sufficient.
\subsection{Improving the bound}
We fix a number field \K. We denote $q_{\K,N,\delta}$ the
restriction of $q_\K$ to $\SC(N,\delta)$. According to
Corollary~\ref{coroKBL}, if $q_{\K,N,\delta}$ has a negative eigenvalue then
$2N\delta$ is a bound for \K. This justifies the following definition.
\begin{defi}
The pair $(N,\delta)$ is \good when $q_{\K,N,\delta}$ has a negative eigenvalue.
\end{defi}

We can reinterpret Functions~\GRHcheck and \BDyDF saying that if
$\GRHcheck(\K,2\delta)$ is negative then $(1,\delta)$ is \good and that
$\left(1,\frac12\log\BDyDF(\K)\right)$ is \good.

As a first step to improve on Corollary~\ref{coroKB}, given $\delta>0$ we look
for the smallest $N$ such that $(N,\delta)$ is \good. Looking for such an $N$
can be done fairly easily with this setup. For any $i>=1$, let $\Phi_i$ be the
characteristic function of $({-i\delta},i\delta)$. Then $(\Phi_i)_{1<=i<=N}$
is a basis of $\SC(N,\delta)$. We have
$\Phi_i\ast\Phi_i=F_{2i\delta}=(2i\delta-|x|)C_{4i\delta}(|x|)$; observe also
that the function considered in Corollary~\ref{coroKB} is
$\frac1{\log T}F_{\log T}$. We further observe that
$$\Phi_i\ast\Phi_j=F_{(i+j)\delta}-F_{|i-j|\delta}\ .$$
This means that the matrix $A_N$ of $q_{\K,N,\delta}$ can be computed by
computing only the values of $\ell_\K(F_{2i\delta})$ for $1<=i<=2N$ and
subtracting those values.

We then stop when the determinant of $A_N$ is negative or when
$2N\delta>=\BDyDF(\K)$. This does not guarantee that we stop as soon as there
is a negative eigenvalue. Indeed, consider the following sequence of
signatures:
$$(0,p,0)\to(1,p,0)\to(1,p,1)\to(0,p+1,2)\to\cdots$$
We should have stopped when the signature was $(1,p,1)$ however the
determinant was zero there. Our algorithm will stop as soon as there is an odd
number of negative eigenvalues (and no zero) or we go above $\BDyDF(\K)$. Such
unfavorable sequence of signatures is however very unlikely and can be ignored
in practice.

The corresponding algorithm is presented in Function~\NDelta. We have added a
limit $N_{\max}$ for $N$ which is not needed right now but will be used later.
In Function~\NDelta, we need to slightly change \GRHcheck to returns the
difference of both sides of Equation~\eqref{homoeq} instead of~\eqref{coroeq}.
Note that $(\Phi_i)$ is a basis adapted to the inclusion~\eqref{incnext} so
that we only need to compute the edges of the matrix $A_N$ at each step. The
test $\det A<0$ in line~\ref{line:detA<0} can be implemented using Cholesky
$LDL^*$ decomposition which is incremental.

One way to use this function is to compute $T=\BDyDF(\K)$ and for some
$N_{\max}>=2$, let $\delta=\frac{\log T}{2N_{\max}}$ and
$N=\NDelta(\K,\delta,N_{\max})$. Using the inclusion~\eqref{incmul}, we see
that $(N,\delta)$ is \good and that $N<=N_{\max}$, so that we have improved
the bound.
\subsection{Adaptive steps}
Unfortunately Function~\NDelta is not very efficient mostly for two reasons. To
explain them and to improve the function we introduce some extra notations.\\
For any $\delta>0$, let $N_\delta$ be the minimal $N$ such that $(N,\delta)$ is
\good. Observe that Function~\NDelta computes $N_\delta$, as long as
$N_\delta<=N_{\max}$ and no zero eigenvalue prevents success. Obviously,
using~\eqref{incnext}, we see that for any $N>=N_\delta$, $(N,\delta)$ is
\good. We have observed numerically that the sequence $N\delta_N$ is roughly
decreasing, i.e. for most values of $N$ we have
$N\delta_N>=(N+1)\delta_{N+1}$.\\
For any $N>=1$, let $\delta_N$ be the infimum of the $\delta$'s such that
$(N,\delta)$ is \good. It is not necessarily true that if $\delta>=\delta_N$
then $(N,\delta)$ is \good, however we have never found a counterexample. The
function $\delta\mapsto\delta N_\delta$ is piecewise linear with
discontinuities at points where $N_\delta$ changes; the function is increasing
in the linear pieces and decreasing at the discontinuities. This means that if
we take $0<\delta_2<\delta_1$ but we have $N_{\delta_2}>N_{\delta_1}$ then we
may have $N_{\delta_2}\delta_2>N_{\delta_1}\delta_1$ so the bound we get for
$\delta_2$ is not necessarily as good as the one for $\delta_1$.\\
The resolution of Function~\NDelta is not very good: going from $N-1$ to $N$
the bound for the norm of the prime ideals is multiplied by $e^{2\delta}$.
This is the first reason reducing the efficiency of the function. The second
one is that if $N_{\max}$ is above $20$ or so, the number
$\delta=\frac{\log\BDyDF(\K)}{2N_{\max}}$ has no specific reason to be near
$\delta_{N_\delta}$; as discussed above, this means that we can get a better
bound for \K by choosing $\delta$ to be just above either $\delta_{N_\delta}$
or $\delta_{1+N_\delta}$. Both reasons derive from the same facts and give a
bound for \K that can be overestimated by at most $2\delta$ for the considered
$N=\NDelta(\K,\delta,N_{\max})$.

To improve the result, we can use once again inclusion~\eqref{incmul} and
determine a good approximation of $\delta_N$ for $N=2^n$. We determine first by
dichotomy a $\delta_0$ such that $(N_0,\delta_0)$ is \good for some $N_0>=1$
(we use $N_0=8$ in our computation). For any $k>=0$, we take $N_{k+1}=2N_k$ and
determine by dichotomy a $\delta_{k+1}$ such that $(N_{k+1},\delta_{k+1})$ is
\good; we already know that $\frac{\delta_k}2$ is an upper bound for
$\delta_{k+1}$ and we can either use $0$ as a lower bound or try to find a
lower bound not too far from the upper bound because the upper bound is
probably not too bad. The algorithm is described in Function~\Bound. It uses
a subfunction $\OptimalT(\K,N,T_l,T_h)$ which returns the smallest integer
$T\in[T_l,T_h]$ such that $\NDelta(\K,\log T/(2N),N)>0$. The algorithm does not
return a bound below those proved in~\label{theo:Phieasynt}
and~\label{coro:Bach4.01}.
\subsection{Further refinements}
To reduce the time used to compute the determinants, we tried to use steps
of width $4\delta$ in $\left[-\frac12\log T,\frac12\log T\right]$ and of
width $2\delta$ in the rest of $\left[-\frac34\log T,\frac34\log T\right]$, to
halve the dimension of $\SC(N,\delta)$. It worked in the sense that we found
substantially the same $T$ faster. However we decided that the total time of
the algorithm is not high enough to justify the increase in code complexity.
\section{Examples}
In this section we will denote $T(\K)$ the result of Function~\BDyDF and
$T_1(\K)$ the result of Function~\Bound.
\subsection{Various fields}
We tested the algorithm on several fields. Let first
$\K=\quot{\QM[x]}{(P)}$ where
$$\catcode`\*=\active\def*{}
P=x^3 + 559752270111028720*x + 55137512477462689.
$$
The polynomial $P$ has been chosen so that for all primes $2<=p<=53$ there
are two prime ideals of norms $p$ and $p^2$. This ensures that there are
lots of small norms of prime ideals.
We have $T(\K)=19162$. There are $2148$ non-zero prime ideals with norms up to
$T(\K)$. We found that $T_1(\K)=11071$ and that there are $1343$ non-zero prime
ideals of norms up to $T_1(\K)$. The time used by Function~\BDyDF was 58ms on
our test computer while the time used by our algorithm was an \emph{additional}
36ms. The test was designed in such a way that our algorithm used the
decomposition information of Function~\BDyDF, so it saved a little time.

\enlargethispage{3\baselineskip}
We tested also the algorithm on the set of $4686$ fields of degree $2$ to $27$
and small discriminant coming from a benchmark of~\cite{PARI}. The mean value
of $\frac{T_1(\K)}{T(\K)}$ for those fields is lower than $\frac12$.

For cyclotomic fields, the new algorithm does not give results significantly
better than those of Belabas, Diaz y Diaz and Friedman. It might be because
the discriminant of a cyclotomic field is not large enough with respect to
its degree.
\subsection{Pure fields}
We computed $T(\K)$ and $T_1(\K)$ for fields of the form $\quot{\QM[x]}{(P)}$
with $P=x^n\pm p$ and $p$ is the first prime after $10^a$ for a certain family
of integers $n$ and $a$. We computed the family of $\frac{T_1(\K)}{T(\K)}$ for
each fixed degree. The graph shows that it is decreasing with the
discriminant. The graph of $\frac{T_1(\K)}{T(\K)}(\llDK)^2$ is much more
regular and looks to have a non-zero limit, see Figure~\ref{fig1} below. We
computed the mean of $\frac{T_1(\K)}{T(\K)}(\llDK)^2$ for each fixed degree.
The results are summarized below:
\[
\let\mc=\multicolumn
\begin{array}{l|r|r|l}
P   & a<={} & \lDK<={}  & \text{mean}   \\
\hline
x^2-p   & 3999  &  9212     & 13.19     \\
x^6+p   & 1199  & 13818     & 13.38     \\
x^{21}-p&  328  & 15169     & 13.68
\end{array}
\]
%G my(a= 2);ldisc=read(Str("disc",a,".ldisc"));buch=readvec(Str("disc",a,".buch"));[#buch,ldisc[#buch],sum(n=1,#buch,buch[n][2,1]/buch[n][1,1]*log(ldisc[n])^2)/#buch]
%G my(a=-6);ldisc=read(Str("disc",a,".ldisc"));buch=readvec(Str("disc",a,".buch"));[#buch,ldisc[#buch],sum(n=1,#buch,buch[n][2,1]/buch[n][1,1]*log(ldisc[n])^2)/#buch]
%G my(a=21);ldisc=read(Str("disc",a,".ldisc"));buch=readvec(Str("disc",a,".buch"));[#buch,ldisc[#buch],sum(n=1,#buch,buch[n][2,1]/buch[n][1,1]*log(ldisc[n])^2)/#buch]
%
The small discriminants are (obviously) much less sensitive to the new
algorithm. We reduced the range for each series to have $\lDK<=500$. The
results are as follows:
\[
\begin{array}{l|r|l}
P   & a<={} & \text{mean}   \\
\hline
x^2-p   &   218 & 12.35     \\
x^6+p   &    43 & 13.66     \\
x^{21}-p&    10 & 17.19
\end{array}
\]
%G my(a= 2);ldisc=read(Str("disc",a,".ldisc"));buch=readvec(Str("disc",a,".buch"));N=vecmin(select(n->ldisc[n]>500,[1..#buch]))-1;[N,vecsum(vector(N,n,(buch[n][2,1]/buch[n][1,1]*log(ldisc[n])^2)))/N]
%G my(a=-6);ldisc=read(Str("disc",a,".ldisc"));buch=readvec(Str("disc",a,".buch"));N=vecmin(select(n->ldisc[n]>500,[1..#buch]))-1;[N,vecsum(vector(N,n,(buch[n][2,1]/buch[n][1,1]*log(ldisc[n])^2)))/N]
%G my(a=21);ldisc=read(Str("disc",a,".ldisc"));buch=readvec(Str("disc",a,".buch"));N=vecmin(select(n->ldisc[n]>500,[1..#buch]))-1;[N,vecsum(vector(N,n,(buch[n][2,1]/buch[n][1,1]*log(ldisc[n])^2)))/N]
%
\subsection{Biquadratic fields}
We repeated the computations above also for biquadratic fields
$\QM[\sqrt{p_1},\sqrt{p_2}]$ where each $p_i$ is the first prime after
$10^{a_i}$ for a certain family of integers $a_i$. We found that the mean of
$\frac{T_1(\K)}{T(\K)}(\llDK)^2$ is $13.63$ for the $7119$ fields computed
and $13.88$ if we restrict the family to the $1537$ ones with $\lDK<=500$.
%
%G my(a="4p");ldisc=read(Str("disc",a,".ldisc"));buch=readvec(Str("disc",a,".buch"));[#buch,ldisc[#buch],sum(n=1,#buch,buch[n][2,1]/buch[n][1,1]*log(ldisc[n])^2)/#buch]
%G my(a="4p");ldisc=read(Str("disc",a,".ldisc"));buch=readvec(Str("disc",a,".buch"));N=vecmin(select(n->ldisc[n]>500,[1..#buch]))-1;[N,vecsum(vector(N,n,buch[n][2,1]/buch[n][1,1]*log(ldisc[n])^2))/N]
%
\subsection*{Final remarks}
In~\cite[Th. 4.3]{small-generators} the authors prove that for a fixed degree
$T(\K)\gg(\lDK\llDK)^2$ and conjecture that
$T(\K)\sim\frac1{16}(\lDK\llDK)^2$ while our computations suggest that
$T_1(\K)$ has smaller order. We will prove in a subsequent
article~\cite{quality} that $T(\K)\asymp(\lDK\llDK)^2$ and that
$T_1(\K)\ll(\log\DK)^2$.
\medskip

\centerline{%
\begin{tabular}{c}
% GNUPLOT: LaTeX picture with Postscript
\begingroup
  \makeatletter
  \providecommand\color[2][]{%
    \GenericError{(gnuplot) \space\space\space\@spaces}{%
      Package color not loaded in conjunction with
      terminal option `colourtext'%
    }{See the gnuplot documentation for explanation.%
    }{Either use 'blacktext' in gnuplot or load the package
      color.sty in LaTeX.}%
    \renewcommand\color[2][]{}%
  }%
  \providecommand\includegraphics[2][]{%
    \GenericError{(gnuplot) \space\space\space\@spaces}{%
      Package graphicx or graphics not loaded%
    }{See the gnuplot documentation for explanation.%
    }{The gnuplot epslatex terminal needs graphicx.sty or graphics.sty.}%
    \renewcommand\includegraphics[2][]{}%
  }%
  \providecommand\rotatebox[2]{#2}%
  \@ifundefined{ifGPcolor}{%
    \newif\ifGPcolor
    \GPcolorfalse
  }{}%
  \@ifundefined{ifGPblacktext}{%
    \newif\ifGPblacktext
    \GPblacktexttrue
  }{}%
  % define a \g@addto@macro without @ in the name:
  \let\gplgaddtomacro\g@addto@macro
  % define empty templates for all commands taking text:
  \gdef\gplbacktext{}%
  \gdef\gplfronttext{}%
  \makeatother
  \ifGPblacktext
    % no textcolor at all
    \def\colorrgb#1{}%
    \def\colorgray#1{}%
  \else
    % gray or color?
    \ifGPcolor
      \def\colorrgb#1{\color[rgb]{#1}}%
      \def\colorgray#1{\color[gray]{#1}}%
      \expandafter\def\csname LTw\endcsname{\color{white}}%
      \expandafter\def\csname LTb\endcsname{\color{black}}%
      \expandafter\def\csname LTa\endcsname{\color{black}}%
      \expandafter\def\csname LT0\endcsname{\color[rgb]{1,0,0}}%
      \expandafter\def\csname LT1\endcsname{\color[rgb]{0,1,0}}%
      \expandafter\def\csname LT2\endcsname{\color[rgb]{0,0,1}}%
      \expandafter\def\csname LT3\endcsname{\color[rgb]{1,0,1}}%
      \expandafter\def\csname LT4\endcsname{\color[rgb]{0,1,1}}%
      \expandafter\def\csname LT5\endcsname{\color[rgb]{1,1,0}}%
      \expandafter\def\csname LT6\endcsname{\color[rgb]{0,0,0}}%
      \expandafter\def\csname LT7\endcsname{\color[rgb]{1,0.3,0}}%
      \expandafter\def\csname LT8\endcsname{\color[rgb]{0.5,0.5,0.5}}%
    \else
      % gray
      \def\colorrgb#1{\color{black}}%
      \def\colorgray#1{\color[gray]{#1}}%
      \expandafter\def\csname LTw\endcsname{\color{white}}%
      \expandafter\def\csname LTb\endcsname{\color{black}}%
      \expandafter\def\csname LTa\endcsname{\color{black}}%
      \expandafter\def\csname LT0\endcsname{\color{black}}%
      \expandafter\def\csname LT1\endcsname{\color{black}}%
      \expandafter\def\csname LT2\endcsname{\color{black}}%
      \expandafter\def\csname LT3\endcsname{\color{black}}%
      \expandafter\def\csname LT4\endcsname{\color{black}}%
      \expandafter\def\csname LT5\endcsname{\color{black}}%
      \expandafter\def\csname LT6\endcsname{\color{black}}%
      \expandafter\def\csname LT7\endcsname{\color{black}}%
      \expandafter\def\csname LT8\endcsname{\color{black}}%
    \fi
  \fi
    \setlength{\unitlength}{0.0500bp}%
    \ifx\gptboxheight\undefined%
      \newlength{\gptboxheight}%
      \newlength{\gptboxwidth}%
      \newsavebox{\gptboxtext}%
    \fi%
    \setlength{\fboxrule}{0.5pt}%
    \setlength{\fboxsep}{1pt}%
\begin{picture}(8502.00,3968.00)%
    \gplgaddtomacro\gplbacktext{%
      \csname LTb\endcsname%
      \put(462,440){\makebox(0,0)[r]{\strut{}$10$}}%
      \put(462,906){\makebox(0,0)[r]{\strut{}$11$}}%
      \put(462,1372){\makebox(0,0)[r]{\strut{}$12$}}%
      \put(462,1838){\makebox(0,0)[r]{\strut{}$13$}}%
      \put(462,2305){\makebox(0,0)[r]{\strut{}$14$}}%
      \put(462,2771){\makebox(0,0)[r]{\strut{}$15$}}%
      \put(462,3237){\makebox(0,0)[r]{\strut{}$16$}}%
      \put(462,3703){\makebox(0,0)[r]{\strut{}$17$}}%
      \put(657,220){\makebox(0,0){\strut{}$0$}}%
      \put(2387,220){\makebox(0,0){\strut{}$4000$}}%
      \put(4117,220){\makebox(0,0){\strut{}$8000$}}%
      \put(5847,220){\makebox(0,0){\strut{}$12000$}}%
      \put(7577,220){\makebox(0,0){\strut{}$16000$}}%
      \put(7709,440){\makebox(0,0)[l]{\strut{}$10$}}%
      \put(7709,906){\makebox(0,0)[l]{\strut{}$11$}}%
      \put(7709,1372){\makebox(0,0)[l]{\strut{}$12$}}%
      \put(7709,1838){\makebox(0,0)[l]{\strut{}$13$}}%
      \put(7709,2305){\makebox(0,0)[l]{\strut{}$14$}}%
      \put(7709,2771){\makebox(0,0)[l]{\strut{}$15$}}%
      \put(7709,3237){\makebox(0,0)[l]{\strut{}$16$}}%
      \put(7709,3703){\makebox(0,0)[l]{\strut{}$17$}}%
    }%
    \gplgaddtomacro\gplfronttext{%
      \colorrgb{0.58,0.00,0.83}%
      \put(5993,3530){\makebox(0,0)[l]{\strut{}quadratic}}%
      \colorrgb{0.00,0.62,0.45}%
      \put(5993,3310){\makebox(0,0)[l]{\strut{}degree 6}}%
      \colorrgb{0.34,0.71,0.91}%
      \put(5993,3090){\makebox(0,0)[l]{\strut{}degree 21}}%
      \colorrgb{0.90,0.62,0.00}%
      \put(5993,2870){\makebox(0,0)[l]{\strut{}biquadratic}}%
    }%
    \gplbacktext
    \put(0,0){\includegraphics{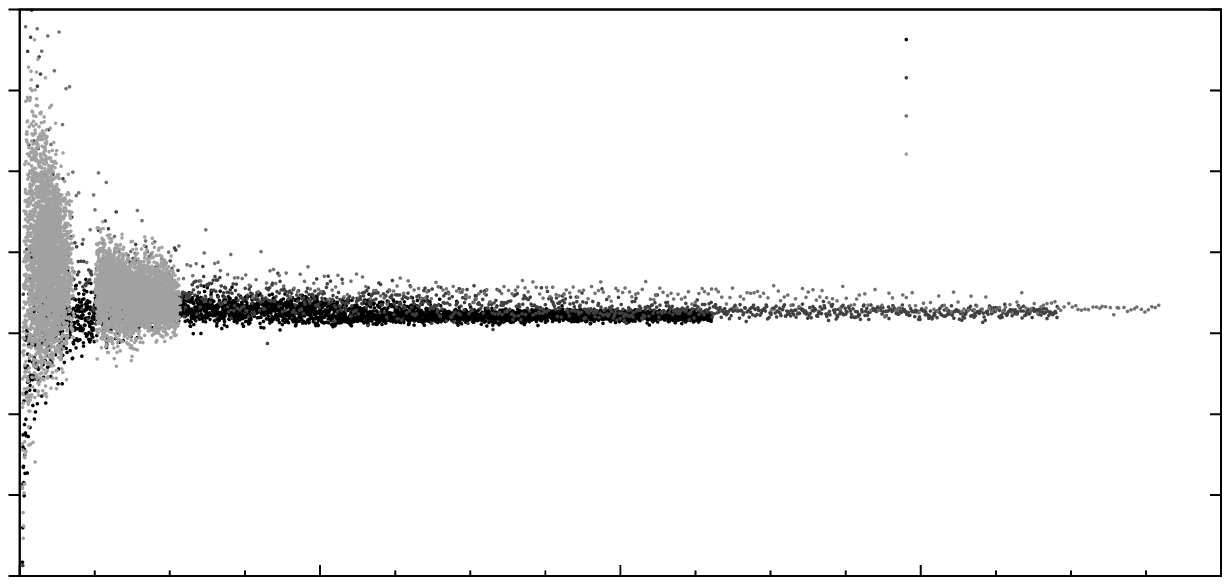}}%
    \gplfronttext
  \end{picture}%
\endgroup
\\
\textsc{Figure 1}: $\frac{T_1(\K)}{T(\K)}(\llDK)^2$ for some pure fields;
in abscissa \lDK.%
{\expandafter\def\csname @currentlabel\endcsname{1}\label{fig1}}
\end{tabular}%
}
\begin{function}
  \SetAlgoLined
  \KwIn{a number field \K}
  \KwIn{a positive real $\delta$}
  \KwIn{a positive integer $N_{\max}$}
  \KwOut{an $N\leqslant N_{\max}$ such that $(N,\delta)$ is \good or $0$}

  $tab\leftarrow\text{$(2N_{\max}+1)$-dimensional array}$\;
  $tab[0]\leftarrow0$\;
  $A\leftarrow\text{$N_{\max}\times N_{\max}$ identity matrix}$\;
  $N\leftarrow 0$\;
  \While{$N<N_{\max}$}{
    $N\leftarrow N+1$\;
    $tab[2N-1]\leftarrow(2N-1)\GRHcheck(\K,(2N-1)\delta)$\;
    $tab[2N]\leftarrow2N\GRHcheck(\K,2N\delta)$\;
    \For{$i\leftarrow1$ \KwTo $N$}{
      $A[N,i]\leftarrow tab[N+i]-tab[N-i]$\;
      $A[i,N]\leftarrow A[N,i]$\;
    }
    \If{$\det A<0$}{\label{line:detA<0}
      \KwRet{$N$}\;
    }
  }
  \KwRet{$0$}\;
  \caption{NDelta(\K,$\delta$,$N_{\max}$)}
\end{function}
\begin{function}
  \SetAlgoLined
  \KwIn{a number field \K}
  \KwOut{a bound for the norm of a system of generators of $\Cl_\K$}
  $T_0\leftarrow 4\left(\lDK+\llDK-(\g+\log 2\pi)\nK+1+(\nK+1)\frac{\log(7\lDK)}\lDK\right)^2$\;
  $T_0\leftarrow \min\left(T_0,4.01\log^2\DK\right)$\;
  $N\leftarrow8$;
  $\delta\leftarrow0.0625$\;
  \While{$\NDelta(\K,\delta,N)=0$}{
    $\delta\leftarrow\delta+0.0625$\;
  }
  $T_h\leftarrow \OptimalT(\K,N,e^{2N\,(\delta-0.0625)},e^{2N\,\delta})$\;
  $T\leftarrow T_h+1$\;
  \While{$T_h<T\mathop{||}T>T_0$}{
    $T\leftarrow T_h$;
    $N\leftarrow2N$\;
    $T_h\leftarrow \OptimalT(\K,N,1,T_h)$\;
  }
  \KwRet{$T$}\;
  \caption{Bound(\K)}
\end{function}

\pagebreak

%\bibliography{mybiblio}
\vfill
{\catcode`\$=10
\def\ident Id: #1,v #2 #3 #4 #5 Exp {\def\thefootnote{}\footnotetext{\tiny version #2, #3 #4}}
\ident$Id: note.tex,v 1.77 2016/09/30 13:48:08 grenie Exp $
}
\end{document}